\documentclass[12pt,reqno]{amsart}
\usepackage{amsmath}
\usepackage{amsthm}
\usepackage{amsfonts}
\usepackage{amssymb}
\usepackage{mathrsfs}
\usepackage[margin=1.35in, top =1.25in, bottom = 1.25in]{geometry}
\usepackage{setspace}
\usepackage{tikz}
\usetikzlibrary{matrix, arrows}
\usepackage{tikz-cd}
\usepackage{enumitem}
\usepackage[colorlinks, breaklinks]{hyperref}
\hypersetup{linkcolor=black, citecolor=blue, filecolor=black, urlcolor=black} 
\usepackage{pdflscape}
\usepackage{array}
\usepackage{adjustbox}
\usepackage{tabularx}
\usepackage{multirow}
\usepackage{multicol}
\usepackage{float}

\usepackage{amsxtra}
\usepackage{amsrefs}
\usepackage[all]{xy}
\usepackage{verbatim}
\usepackage{hyphenat}

\DeclareMathOperator{\ch}{char}
\newcommand{\FF}{\mathbb{F}}			
\DeclareMathOperator{\GL}{GL}			

\newcommand{\KK}{\mathbb{K}}			

\DeclareMathOperator{\SL}{SL}			
\renewcommand{\vec}[1]{\mathbf{#1}}		
\newcommand{\ZZ}{\mathbb{Z}}

\theoremstyle{plain}

\theoremstyle{definition}

\newtheorem*{ack}{Acknowledgments}

\theoremstyle{remark}

\title{The InvariantRing package for Macaulay2}
\author[Ferraro, Galetto, Gandini, Huang, Mastroeni, and Ni]{Luigi Ferraro, Federico Galetto, Francesca Gandini, Hang Huang, Matthew Mastroeni, and Xianglong Ni}
\date{}

\begin{document}

\maketitle

\begin{abstract}
  We describe a significant update to the existing
  \texttt{InvariantRing} package for \textit{Macaulay2}.  In addition
  to expanding and improving the methods of the existing package for
  actions of finite groups, the updated package adds functionality for
  computing invariants of diagonal actions of tori and finite abelian
  groups as well as invariants of arbitrary linearly reductive group
  actions.  The implementation of the package has been completely
  overhauled with the aim of serving as a unified resource for
  invariant theory computations in \textit{Macaulay2}.
\end{abstract}

\section{Introduction}
Let $G$ be a group acting linearly on an $n$-dimensional vector space $V$ over a field $\KK$ via $\vec{v} \mapsto g \cdot \vec{v}$ for $g \in G$ and $\vec{v} \in V$. The action of $G$ on $V$ induces an action of $G$ on the polynomial ring $\KK[V] = \KK[x_1, \ldots, x_n]$ by $(g \cdot f)(\vec{x}) = f(g^{-1} \cdot \vec{x})$. By a classical result of Hilbert (see \cite[2.2.10]{Derksen:Kemper}), the subring of invariant polynomials
\[ \KK[V]^G = \{ f \in \KK[V] \mid g \cdot f = f, \, \forall g \in G \} \]
 is always a finitely generated $\KK$-algebra provided that $G$ is linearly reductive and $V$ is a rational representation of $G$.\footnote{However, this is not the case in general by Nagata's celebrated solution to Hilbert's Fourteenth Problem \cite{Nagata}.}  A linearly reductive group is a group that can be identified with a Zariski-closed subgroup of some general linear group $\GL_n = \GL_n(\KK)$ and has good representation-theoretic properties while still being general enough to encompass all finite groups, all tori $(\KK^\times)^r$, and all semisimple Lie groups (at least if $\ch(\KK) = 0$).  Since Hilbert's proof of the finite generation of such invariant rings is famously nonconstructive, determining a minimal set of generating invariants for a specific linearly reductive group action remains a central question in invariant theory.  This article describes a significant update\footnote{\url{https://github.com/galettof/InvariantRing}} to the package \texttt{InvariantRing} for \textit{Macaulay2} \cite{Macaulay2}, originally created by Hawes \cite{Hawes}, implementing several algorithms to compute generators for $\KK[V]^G$ for various types of group actions. 
 
 \begin{ack}
 Much of the work on the updated version of this package was completed as part of the virtual \textit{Macaulay2} workshop in May 2020, originally to be held at Cleveland State University, which was partially supported by NSF award DMS-2003883.  We thank Harm Derksen and Visu Makam for their helpful suggestions and feedback on the package.
 \end{ack}

\section{Types of Group Actions}
\label{sec:functionality}

The updated \texttt{InvariantRing} package introduces the \texttt{GroupAction} type and the specialized subtypes \texttt{FiniteGroupAction}, \texttt{DiagonalAction}, and \texttt{Linearly}\texttt{Reduct}\-\texttt{iveAction} as a unified approach to creating group actions in \textit{Macaulay2}. Each type has its own constructor and methods for computing invariants.


\subsection{Finite Group Actions}
\label{sec:finite-groups-actions}

Creating a group action of type \texttt{FiniteGroupAction} requires a polynomial ring $R = \KK[x_1, \dots, x_n]$ over a field and a list of $n \times n$ matrices generating a finite subgroup of $\GL_n$. The example below illustrates the construction of the natural action of the alternating group $A_4$ on a polynomial ring in four variables.

%
%

%

\begin{verbatim}
i2 : R = QQ[x_1..x_4];

i3 : L = {permutationMatrix "2314", permutationMatrix "2143"};

i4 : A4 = finiteAction(L, R)

o4 = R <- {| 0 0 1 0 |, | 0 1 0 0 |}
           | 1 0 0 0 |  | 1 0 0 0 |
           | 0 1 0 0 |  | 0 0 0 1 |
           | 0 0 0 1 |  | 0 0 1 0 |

o4 : FiniteGroupAction

\end{verbatim}

A minimal set of algebra generators for the ring of invariants can be computed using the method \texttt{invariants}. 
\begin{verbatim}

i5 : invariants A4

                          2    2    2    2   3    3    3    3   4
o5 = {x  + x  + x  + x , x  + x  + x  + x , x  + x  + x  + x , x 
       1    2    3    4   1    2    3    4   1    2    3    4   1
     -------------------------------------------------------------
        4    4    4   3 2        3 2    2   3    2 3      3 2    
     + x  + x  + x , x x x  + x x x  + x x x  + x x x  + x x x  +
        2    3    4   1 2 3    1 2 3    1 2 3    1 2 4    1 3 4  
     -------------------------------------------------------------
      2 3      3   2    3   2      3 2      2 3    2   3      2 3
     x x x  + x x x  + x x x  + x x x  + x x x  + x x x  + x x x }
      2 3 4    1 2 4    2 3 4    1 3 4    1 2 4    1 3 4    2 3 4

o5 : List

\end{verbatim}
By default, the invariants of a finite group action are now computed using an implementation of King's algorithm \cite[3.8.2]{Derksen:Kemper} based on the Reynolds operator.  
An alternative algorithm based on linear algebra \cite[\S 3.1.1]{Derksen:Kemper}, which can be faster for large groups, is also available through the option \texttt{UseLinearAlgebra}.

The methods \texttt{molienSeries}, \texttt{primaryInvariants}, and \texttt{secondaryInvariants} from version 1.1.0 of the package are still available. The method previously known as \texttt{invariants}, which computes primary and secondary invariants at once, has been renamed \texttt{hironakaDecomposition}.
We note that a generating set for the ring of invariants can be obtained from a list of primary and secondary invariants, although this is typically more time consuming. For example, on the same machine, the previous computation of invariants for $A_4$ took 0.49 seconds using King's algorithm, 0.05 seconds using linear algebra, and 16.49 seconds by computing a Hironaka decomposition.

\subsection{Diagonal Actions}
\label{sec:diagonal-actions}

A diagonal action over an algebraically closed field $\KK$ is an action by a group $G = T \times A$ that is the product of a torus $T$ and a finite abelian group $A$.  If $T = (\KK^\times)^r$ and $A = \ZZ/d_1\ZZ \times \cdots \times \ZZ/d_s\ZZ$, the action of $G$ on the polynomial ring $R = \KK[x_1, \dots, x_n]$ is \emph{diagonal} if there is an $n \times (r + s)$ matrix of integers $W = (w_{i,j})$, called the \emph{weight matrix}, such that
\[ t \cdot x_j = t_1^{w_{1,j}} \cdots t_r^{w_{r,j}}x_j \]
for all $j$ and all $t = (t_1, \dots, t_r) \in T$, and there is a primitive  $d_i$-th root of unity $\zeta_i$  such that the generator $u_i$ of the cyclic abelian factor $\ZZ/d_i\ZZ$ acts by
\[ u_i \cdot x_j = \zeta_i^{w_{r+i,j}} x_j \] 
for all $i, j$.  Because a diagonal action preserves the natural $\ZZ^n$-grading of $R$, a polynomial is easily seen to be invariant under a diagonal action if and only if each of its terms is invariants.  Thus, we can choose a set of invariant monomials as a minimal set of generating invariants. 

Creating a group action of type \texttt{DiagonalAction} requires a polynomial ring $R$, a weight matrix $W$, and list of positive integers defining the orders of the cyclic factors of the group $G$.  We illustrate this below for a torus action over $\FF_9$. \vspace{1 ex}
\begin{verbatim}
i2 : R = (GF 9)[x_1..x_4];

i3 : W = matrix {{5, -3, -1, 4}, {-3, 1, 1, 5}, {0, -4, 2, 6}};

              3        4
o3 : Matrix ZZ  <--- ZZ

i4 : T = diagonalAction(W, {}, R)

                 * 3
o4 = R <- ((GF 9) )  via 

     | 5  -3 -1 4 |
     | -3 1  1  5 |
     | 0  -4 2  6 |

o4 : DiagonalAction

i5 : invariants T

           2
o5 = {x x x }
       1 2 3

o5 : List
\end{verbatim}

Because a diagonal action involving finite cyclic factors requires the existence of roots of unity for the action to be defined, the invariants computed for a diagonal action should \emph{always} be understood to be the minimal generating invariant monomials over an appropriate infinite extension of the coefficient field.  Over such an extension field, the monomial $\vec{x}^\vec{a} = x_1^{a_1} \cdots x_n^{a_n}$ is invariant if and only if the $i$-th coordinate of the vector $W\vec{a}$ is zero for $1 \leq i \leq r$ and the $(r + i)$-th coordinate is zero modulo $d_i$ for $1 \leq i \leq s$.  As a result, finding a minimal set of invariant monomials for a diagonal action reduces to a suitable problem in polyhedral geometry.  Using the fact that $R^G = (R^A)^T$, we have implemented a recent algorithm of Gandini \cite{Gandini} for first computing invariants of the finite abelian part of a diagonal action; we then use a modified version of \cite[4.3.1]{Derksen:Kemper} to find the torus-invariant monomials in $R^A$.

As noted above, some form of polyhedral geometry computation is unavoidable when computing diagonal action invariants.  The current implementation relies on one of the \textit{Macaulay2} packages \texttt{Polyhedra} \cite{Polyhedra} or \texttt{Normaliz} \cite{Normaliz:M2} to handle this part of the computation.  The latter package is simply an interface to the stand-alone Normaliz program \cite{Normaliz} written in C++.  As a program specializing in polyhedral geometry, Normaliz already has the ability to compute the invariants of a diagonal action.  Although Normaliz can have a significant speed advantage on larger examples since it is written in C++ -- for example, the above computation on the same machine takes about 25.60 seconds using the Derksen-Gandini algorithm versus 0.01 seconds by calling Normaliz directly -- we believe there is value in having an algorithm for computing diagonal invariants implemented in \textit{Macaulay2} that is independent of any particular external software.

In addition, since \textit{Macaulay2} has the ability to perform computations over any finite field, it is possible to define diagonal actions (as in the above example) where the action makes sense over the coefficient field specified by the user rather than a suitable infinite extension.  In this case, our algorithm includes the ability to compute a minimal set of generating monomials \emph{literally} over the given coefficient field through the option \texttt{UseCoefficientRing}, a feature that Normaliz currently lacks.  Computing invariants literally over $\FF_9$ in the above example, we find the additional invariants:
\vspace{1 em}
\begin{verbatim}

i6 : invariants(T, UseCoefficientRing => true)

           2   8   8   4 4   8   2 6   4 4   4 4   6 2   8
o6 = {x x x , x , x , x x , x , x x , x x , x x , x x , x }
       1 2 3   4   3   2 3   2   1 2   1 3   1 2   1 2   1

o6 : List

\end{verbatim}
We note that the above computation took only 0.53 seconds on the same machine as the previous calculation.

\subsection{Linearly Reductive Actions}
\label{sec:line-reduct-actions}

In the case of a group action by an arbitrary linearly reductive group $G$, an action of type \texttt{LinearlyReductiveAction} is constructed from the data of an ideal $I$ in an ambient polynomial ring $S$ defining the linearly reductive group $G$ as a Zariski-closed subset of some $\KK^m$, a polynomial ring $R = \KK[x_1, \dots, x_n]$ on which the group acts, and an $n \times n$ matrix of polynomials in $S$ defining the action of $G$ on $R$.

The example below illustrates how to set up the classic action of $\SL_2$ on the coefficients of binary quadrics $ax^2 + bxy + cy^2$ in the ring $\KK[x, y]$ where $\ch(\KK) = 0$. We begin by constructing the defining data for the group $\SL_2$.
\vspace{1 em}
\begin{verbatim}
i2 : S = QQ[z_(1,1)..z_(2,2)];

i3 : I = ideal(z_(1,1)*z_(2,2) - z_(1,2)*z_(2,1) - 1)

o3 = ideal(- z   z    + z   z    - 1)
              1,2 2,1    1,1 2,2

o3 : Ideal of S

\end{verbatim}
There is a natural action of $\SL_2$ on $\KK[x, y]$ by linear changes of coordinates.  We construct the matrix representing the restriction of this action to the space of quadrics $\KK[x, y]_2$ with respect to the monomial basis $x^2, xy, y^2$. \vspace{1 em}
\begin{verbatim}
i4 : A = S[x,y];

i5 : M = (map(S,A)) last coefficients sub(basis(2,A),
         {x => z_(1,1)*x+z_(1,2)*y, y => z_(2,1)*x+z_(2,2)*y});

             3       3
o5 : Matrix S  <--- S

\end{verbatim}
Viewing $R = \KK[a, b, c]$ as the dual polynomial ring of coefficients of the quadrics in $\KK[x, y]_2$, the transpose of this matrix represents the induced $\SL_2$-action on $R$. \vspace{1 em}
\begin{verbatim}
i6 : R = QQ[a,b,c];

i7 : L = linearlyReductiveAction(I, transpose M, R)

o7 = R <- S/ideal(- z   z    + z   z    - 1) via 
                     1,2 2,1    1,1 2,2

     | z_(1,1)^2      2z_(1,1)z_(1,2)               z_(1,2)^2      |
     | z_(1,1)z_(2,1) z_(1,2)z_(2,1)+z_(1,1)z_(2,2) z_(1,2)z_(2,2) |
     | z_(2,1)^2      2z_(2,1)z_(2,2)               z_(2,2)^2      |

o7 : LinearlyReductiveAction

\end{verbatim}

The \emph{Hilbert ideal} of the action is the ideal of the polynomial ring $R$ generated by all non-constant homogeneous invariant polynomials.  A minimal set of invariant generators for the Hilbert ideal forms a minimal set of algebra generators for the ring of invariants.  The package uses an implementation of \cite[4.2.8]{Derksen:Kemper} to compute a set of generators of the Hilbert ideal.
\vspace{1 em}
\begin{verbatim}
i8 : hilbertIdeal L

            2
o8 = ideal(b  - 4a*c)

o8 : Ideal of R

\end{verbatim}
In this case, the generator is none other than the discriminant, and it is already invariant.  Typically, this is not the case, and so, the final step is to find \emph{invariant} generators of the Hilbert ideal.  As the Reynolds operator is not yet implemented for a general linearly reductive group (see last section), the package uses an alternative approach \cite[4.5.1]{Derksen:Kemper} to find a vector space basis of the invariants degree-by-degree.

\section{Rings of Invariants}
\label{sec:common-methods}

Finally, the package also introduces the \texttt{RingOfInvariants} type as a container for all ring-theoretic information about a given group action.  Rings of invariants can be computed by calling the method \texttt{invariantRing} on any type of group action or by using the natural superscript notation.  We illustrate the latter for the finite group action of $A_4$ on a polynomial ring $R$ defined in \S \ref{sec:finite-groups-actions}.
\vspace{1 em}
\begin{verbatim}
i6 : S=R^A4

o6 =                         2    2    2    2   3    3    3    3 
      QQ[x  + x  + x  + x , x  + x  + x  + x , x  + x  + x  + x ,
          1    2    3    4   1    2    3    4   1    2    3    4 
      -----------------------------------------------------------
       4    4    4    4   3 2        3 2    2   3    2 3    
      x  + x  + x  + x , x x x  + x x x  + x x x  + x x x  +
       1    2    3    4   1 2 3    1 2 3    1 2 3    1 2 4  
      -----------------------------------------------------------
       3 2      2 3      3   2    3   2      3 2      2 3  
      x x x  + x x x  + x x x  + x x x  + x x x  + x x x  +
       1 3 4    2 3 4    1 2 4    2 3 4    1 3 4    1 2 4  
      -----------------------------------------------------------
       2   3      2 3
      x x x  + x x x ]
       1 3 4    2 3 4
       
o6 : RingOfInvariants

\end{verbatim}

The \texttt{RingOfInvariants} type has access to methods such as \texttt{definingIdeal}, which gives a presentation of the ring of invariants as an affine algebra, and \texttt{hilbertSeries}.  Unlike the Hilbert series of a typical quotient ring, the Hilbert series of a ring of invariants of a finite group action is presented in a way that emphasizes the degrees of the primary invariants, which generate a polynomial subring of the ring on invariants over which it is finitely generated as a module.

\begin{verbatim}
i7 : hilbertSeries(S,Reduce=>true)

                       6
                  1 + T
o7 = -------------------------------
           4       3       2
     (1 - T )(1 - T )(1 - T )(1 - T)

o7 : Expression of class Divide
\end{verbatim}

\section{Future Directions}
The major drawback of the function \texttt{hilbertIdeal} is that, even though its output generates all invariants over the polynomial ring, the generators themselves do not need to be invariant. A set of algebra generators for the ring of invariants can then be obtained by applying the Reynolds operator of a linearly reductive group to the generators of the Hilbert ideal.  Currently, the method \texttt{reynoldsOperator} implements the Reynolds operator for finite group and diagonal actions.  However, the Reynolds operator has not yet been implemented for infinite linearly reductive groups because its definition requires building the Lie algebra structure associated to the group. One special exception where it is possible to explicitly construct the Reynolds operator is the case of $\GL_n$ acting on a vector space $V$.  In this case, Cayley's Omega Process \cite[\S 4.5.3]{Derksen:Kemper} constructs the Reynolds operator in terms of the multiplication map on coordinate rings $\mu^* : \KK[V] \to \KK[V] \otimes_\KK \KK[\GL_n]$ and all partial derivatives on the coordinate ring $\KK[\GL_n]$. We plan to implement Cayley's Omega Process in a future update. More generally, implementing Reynolds operators for all linearly reductive groups is a long-term goal for the development of the \texttt{InvariantRing} package.

\end{document}